\documentclass[12pt]{article}

\usepackage{amsmath,amssymb,amsfonts,theorem,makeidx,latexsym,epsfig,,subfigure}

\newtheorem{defn}{Definition}[section]

\newtheorem{lemma}[defn]{Lemma}

{\theorembodyfont{\rmfamily}

	\newtheorem{ex}[defn]{Example}}

\newtheorem{thm}[defn]{Theorem}

\newtheorem{prop}[defn]{Proposition}

\newtheorem{cor}[defn]{Corollary}

\numberwithin{equation}{section}

\newcommand{\h}{{\cal H}}

\newcommand{\ltr}{ L^2(\mathbb R) }

\newcommand{\mn}{\mathbb N}

\newcommand{\mr}{\mathbb R}

\newcommand{\mz}{\mathbb Z}

\newcommand{\mc}{\mathbb C}

\def\bp{{\noindent\bf Proof. \ }}

\def\ep{\hfill$\square$\par\bigskip}

\def\bqs{\begin{equation}}

	\def\eqs{\tag*{$\square$}\end{equation}\par\bigskip}

\def\la{\langle}

\def\ra{\rangle}

\def\ftk{\{f_k\}_{k=1}^\infty}

\def\ctk{\{c_k\}_{k=1}^\infty}

\def\gtk{\{g_k \}_{k=1}^\infty}

\def\htk{\{h_k\}_{k=1}^\infty}

\def\etk{\{e_k\}_{k=1}^\infty}

\def\dtk{\{\delta_k\}_{k=1}^\infty}

\def\suk{\sum_{k=1}^\infty}

\def\sukz{\sum_{k\in \mz}}

\def\nl{\left|\left|}

\def\nr{\right|\right|}

\def\span{\overline{\text{span}}}

\def\Span{\text{span}}

\def\supp{\text{supp}}

\def\bop{\begin{op}\rm}

	\def\eop{\end{op}}

\def\bee{\begin{eqnarray}}

	\def\ene{\end{eqnarray}}

\def\bes{\begin{eqnarray*}}

	\def\ens{\end{eqnarray*}}

\def\bei{\begin{itemize}}

	\def\eni{\end{itemize}}

\def\bt{\begin{thm}}

	\def\et{\end{thm}}

\def\bc{\begin{cor}}

	\def\ec{\end{cor}}

\def\bpr{\begin{prop}}

	\def\epr{\end{prop}}

\def\bl{\begin{lemma}}

	\def\el{\end{lemma}}

\def\bd{\begin{defn}}

	\def\ed{\end{defn}}

\def\bex{\begin{ex}}

	\def\enx{\end{ex}}

\def\bfi{\begin{fig}}

	\def\efi{\end{fig}}

\def\sukz{\sum_{k\in \mz}}

\title{On approximate operator representations of sequences in Banach spaces}

\date{\today}

\author{Ole Christensen\footnote{ 
		Technical University of Denmark,
		DTU Compute,
		Building 303, 2800 Lyngby,
		Denmark,
		Email: ochr@dtu.dk }, Marzieh Hasannasab\footnote{
		Institute of Mathematics,
		TU Berlin,
		Straße des 17. Juni 136, 
		D-10623 Berlin, Germany,
		Email: hasannas@math.tu-berlin.de}, Gabriele Steidl\footnote{
		Institute of Mathematics,
		TU Berlin,
		Straße des 17. Juni 136, 
		D-10623 Berlin, Germany,
		Email: Steidl@math.tu-berlin.de}}

\begin{document}
	
	\maketitle

%
%
%
%


\begin{abstract}
Generalizing results by Halperin et al., Grivaux recently
showed that any linearly independent sequence $\ftk$ in a separable Banach space $X$ can
be represented as a suborbit $\{T^{\alpha(k)}\varphi\}_{k=1}^\infty$ of
some bounded operator $T: X\to X.$ In general, the operator $T$ and
the powers $\alpha(k)$ are not known explicitly. In this paper we consider
approximate representations $\ftk \approx \{T^{\alpha(k)}\varphi\}_{k=1}^\infty$ of certain types of sequences $\ftk;$ in contrast to the results in the literature
we are able to be very explicit about the operator $T$ and suitable powers
$\alpha(k),$ and
we do not need to assume that the sequences are  linearly independent. The
exact meaning of approximation is defined in a way such that $\{T^{\alpha(k)}\varphi\}_{k=1}^\infty$ keeps essential features of $\ftk,$ e.g.,
in the setting of atomic decompositions and Banach frames.
We will present two different approaches. The first approach is universal,
in the sense that it applies in general
Banach spaces;  the technical conditions are typically easy to verify in
sequence spaces, but are more complicated in function spaces. For
this reason we present a second approach, directly tailored to the
setting of Banach function spaces. A number of examples prove that
the results apply in arbitrary weighted $\ell^p$-spaces and $L^p$-spaces.\\


 \end{abstract}
 
\maketitle

\section{Introduction}

A classical result by Halperin, Kitai, and Rosenthal \cite{kitai} states that
if $\ftk$ is any linearly independent sequence in a separable Hilbert space
$\h,$ then there exists a bounded linear operator $T: \h \to \h$ and
appropriate choices of powers $\alpha(k) \in \mn_0$ such that
$f_k= T^{\alpha(k)}f_1,$ for all $k\in \mn.$ This fundamental
observation was later generalized to Banach spaces by Grivaux \cite{Grivaux},
using a different technique. Neither the
operator $T$ nor the appropriate powers $\alpha(k)$ are given in an easily accessible
form in \cite{Grivaux,kitai}.  Generalizing ideas presented in \cite{olemarzieh-10}
in the setting of frames in Hilbert spaces, we will provide an alternative
approach to the questions by Halpering et al. by considering {\it approximate
	operator representations} in a given Banach space $X.$  In other words, we will give up the requirement
that the operator $T$ leads to an exact representation of the given sequence
$\ftk.$ Instead, we will aim at a construction of a bounded operator $T,$
a vector $\varphi \in X,$
and appropriate powers $\alpha(k)$ such that the sequence
$\{T^{\alpha(k)} \varphi\}_{k=1}^\infty$ approximates the given sequence $\ftk$
in various senses to be specified below. We will show that in several  cases we
can specify as well the operator $T$, the vector $\varphi,$ as the powers $\alpha(k).$

The paper is organized as follows. In the rest of the introduction we set the stage by providing
basic definitions and results concerning Banach sequence spaces and
shift-operators on Banach spaces. In Sections \ref{181205a}--\ref{181205b} we
provide general results for obtaining approximate representations of certain
sequences, in the setting of general Banach spaces having a basis. The results are based on
the assumption that the left/right-shift operators with respect to a certain
basis -- see \eqref{050519b} and \eqref{050519ba} below -- are bounded. This
condition can easily be checked in several types of sequence spaces; in particular,
given any weighted $\ell^p$-space, we can specify a scaling of the canonical unit vectors
that satisfy the condition.  We
can also verify the condition in certain Banach spaces of functions; however,  the condition is
difficult to handle in general function spaces. For this reason we provide
an alternative approach, tailored to Banach function spaces, in Section
\ref{190603a}. Here, the assumption of the left/right-shift operators being
bounded is replaced by the condition that the translation operators act
boundedly on the given function space, a condition that is trivially
satisfied in, e.g., all weighted $L^p$-spaces with an
$m$-moderate weight function.
In Section \ref{181205c} we show
how to construct so-called $\epsilon$-close approximations $\{T^{\alpha(k)} \varphi\}_{k=1}^\infty$ of the sequences $\ftk$ discussed in Sections \ref{181205a}--\ref{190603a}; this paves the way for the results in
Section \ref{181205d}, where it is shown how to construct
the sequence $\{T^{\alpha(k)} \varphi\}_{k=1}^\infty$ such that it keeps key
features of the sequence $\ftk$ in the setting of atomic decompositions.

In the entire paper $X$ will denote a separable Banach space, and   $X_d$ will be a Banach
space consisting of scalar-valued sequences indexed by $\mn.$ We will refer to such
a space $X_d$ as a {\it Banach sequence space.} We will need the following standard concept
related to Banach sequence spaces.

\bd Let $X_d$ denote a Banach sequence space.

\bei
\item[(i)] $X_d$ is said to be solid if whenever $\ctk \in X_d$ and
$\{b_k\}_{k=1}^\infty$ is any scalar-valued sequence such that
$|b_k| \le |c_k|$ for all $k\in \mn,$ it follows that $\{b_k\}_{k=1}^\infty\in X_d$
and $|| \{b_k\}_{k=1}^\infty|| \le || \ctk ||.$
\item[(ii)] $X_d$ is said to have an absolutely continuous norm if
$|| \{c_k - c_k \chi_{I_n}(k)\}_{k=1}^\infty||\to 0$ as $n\to \infty$
for any sequence $\ctk \in X_d$ and any family of subsets $I_n \subset \mn$ such
that $I_1 \subset I_2 \dots \uparrow \mn.$
\eni
\ed

The above concepts have parallel versions in Banach function space, defined
by obvious modifications; we will apply these without further comments in
the sequel.

\subsection{Shift operators on Banach spaces} \label{181205g}
Recall that a sequence $\etk$ in $X$ is  a (Schauder) {\it basis} if every element $x\in X$ has a unique representation  $x=\suk c_k e_k$ for some $c_k\in \mc.$
It was proved by Enflo  \cite{enflo} that not all the separable Banach spaces have a basis.
Let $X$ denote a Banach space having a basis $\etk$  and
consider the left-shift and right-shift operators
$L,R: \Span \etk \to X$ given by
\bee\label{050519b}
Le_1=0,\quad L e_k=  e_{k-1}, \quad  k\geq 2, \ene
respectively,
\bee \label{050519ba}
Re_k= e_{k+1},\quad k\in\mn.
\ene
Throughout the paper we will need that the operators $L,R$ extend to bounded linear operators on $X$. This condition in general depends not only on the space $X$ but also on the choice of the basis $\etk$. We will show that the condition is satisfied if $\etk$ is a $p$-Riesz basis,
a concept to be defined next.

\bd Fix some $p\in [1, \infty).$ A sequence $\etk \subset X$ is called a
	$p$-Riesz basis for $X$ if $\span \etk = X$ and there exist constants $A,B>0$ such that
	\bee \label{190609a} A \left(\sum_{k=1}^N |c_k|^p \right)^{1/p} \le \left| \left| \sum_{k=1}^N c_k e_k \right| \right|
	\le B \left(\sum_{k=1}^N |c_k|^p \right)^{1/p},\ene for all finite scalar sequences $\{c_k\}_{k=1}^N$, $N\in\mn$.
	The numbers $A,B$ are called lower, resp. upper bounds.\ed

Typically the lower condition in \eqref{190609a} is more involved to verify than
the upper condition. A convenient criteria for the $p$-Riesz basis property,
which avoids worrying about the lower bound, is stated next.

\bl Assume that $\etk$ is a basis for a reflexive Banach space $X,$ with
dual basis $\{e_k^*\}_{k=1}^\infty.$ Assume that for some $p,q>1$ with $p^{-1}+q^{-1}=1$
there exists a constant $B>0$ such that
\bes \left( \suk |\la f, e_k^* \ra|^p\right)^{1/p} \le B\, ||f||, \,\quad \forall f\in X\ens
and
\bes \left( \suk |\la g, e_k \ra|^q \right)^{1/q} \le B\, ||g||, \, \quad \forall  g\in X^*.\ens
Then $\etk$ is a $p$-Riesz basis for $X$ and $\{e_k^*\}_{k=1}^\infty$ is a $q$-Riesz basis for $X^*$. \el

\bp Let $N\in\mn$ and $\{c_k\}_{k=1}^N$ be a finite scalar sequence. It follows  by H\"{o}lder's inequality that
\bes \| \sum_{k=1}^N c_k e_k \| &=& \sup_{g\in X^*, \| g\|=1} \Big| \la \sum_{k=1}^N c_k e_k , g\ra\Big| \leq \sup_{g\in X^*, \| g\|=1}  \sum_{k=1}^N |c_k| \, \Big|\la  e_k , g\ra \Big| \\
&\leq & \sup_{g\in X^*, \| g\|=1}  \left(\sum_{k=1}^N |c_k|^p \right)^{1/p}  \left(\sum_{k=1}^N |\la g, e_k \ra|^q \right)^{1/q}\\
&\leq & B  \left(\sum_{k=1}^N |c_k|^p \right)^{1/p}.
\ens
Also letting $f:=\sum_{k=1}^N c_k e_k$, we have that $c_k=\la f , e_k^* \ra$ for $k=1,\dots,n$. Therefore
\bes  \left(\sum_{k=1}^N |c_k|^p \right)^{1/p} = \left(\sum_{k=1}^N |\la f , e_k^* \ra|^p \right)^{1/p} \leq B\,\| f\|
= B\,\| \sum_{k=1}^N c_k e_k \|.
\ens
This proves that $\etk$ is a $p$-Riesz basis for $X$ with bounds $B^{-1}$ and $B$. The proof  that $\{e_k^*\}_{k=1}^\infty$ is a $q$-Riesz basis for $X^*$ is similar.
\ep

We now prove that if $\etk$ is a $p$-Riesz basis for some $p\in [1, \infty),$
then indeed the left/right-shift operators with respect to $\etk$ are bounded:

\bpr \label{190609b}  If $\etk$ is a $p$-Riesz basis for $X$
with bounds $A,B$ for some $p\in [1, \infty),$  then the
operators $L,R$ in \eqref{050519b}
and \eqref{050519ba}   extend to bounded linear operators on $X$, and
\bes || L|| \le \frac{B}{A} \, \hspace{1cm}  \, \mbox{and} \, \hspace{1cm} \,
\frac{A}{B} \le || R|| \le \frac{B}{A}.\ens \epr
\bp Given any finite sequence $\{c_k\}_{k=1}^N$,
\bes  \left| \left| L \sum_{k=1}^N c_k e_k \right| \right| & = &
\left| \left| \sum_{k=2}^{N} c_k e_{k-1} \right| \right|
=  \left| \left| \sum_{k=1}^{N-1} c_{k+1} e_{k} \right| \right| \\ & \le &
B\left( \sum_{k=1}^{N-1} |c_{k+1}|^p \right)^{1/p}
\le
B\left( \sum_{k=1}^{N} |c_{k}|^p \right)^{1/p} \\ & \le & \frac{B}{A}
\left| \left|\sum_{k=1}^N c_k e_k \right| \right|.
\ens Thus $L$ extends to a bounded operator on $X,$ with the claimed
estimate of the norm. The proof for the boundedness of the  right-shift operator
and the upper estimate on its norm is similar. The lower bound on
the norm of the operator $R$ follows from $L$ being a left-inverse of $R,$ i.e., $||f||= ||LRf|| \le || L|| \, ||Rf|| \le BA^{-1}||Rf|| \, $ for all $f\in X.$ \ep

As an application of Proposition \ref{190609b} we will now consider
weighted $\ell^p$-spaces. Fixing any $p\in [1,\infty)$ and considering a sequence of positive scalars $\{w_k\}_{k=1}^\infty$, define the  space $\ell^p_w$ by
\bes
\ell^p_w:=\left\{ \ctk \, \big| \, c_k\in\mc\text{ and }\suk |c_k|^p w_k <\infty \right\}.
\ens
Obviously, $\ell_w^p$ is a Banach space with respect to the norm
\bes ||\ctk ||_{p, w}:= \left( \suk  |c_k|^p w_k\right)^{1/p}.\ens
It is also clear that the canonical unit vectors $\{\delta_k\}_{k=1}^\infty$
form a basis for $\ell_w^p.$ The following result shows that in any weighted $\ell^p$-space, we can
specify a certain scaling of the canonical unit vectors that makes
the left/right-shift operators bounded.

\bc \label{190610a} Fix any $ p\in[1,\infty)$, consider a sequence of positive scalars $\{w_k\}_{k=1}^\infty,$ and let $\{\delta_k\}_{k=1}^\infty$
denote the canonical unit basis for $\ell^p_w$. Let $e_k:= \omega_k^{-1/p}\, \delta_k.$
Then $\etk$ is a $p$-Riesz basis for $\ell_w^p$ with bounds $A=B=1.$ In particular,
the left/right-shift operators $L,R$ with respect to the
basis $\etk$ are bounded, and $||L||=||R||=1.$\ec

\bp A simple calculation shows that for any $N\in \mn$ and any finite
scalar sequence $\{c_k\}_{k=1}^N,$
\bes\nl \sum_{k=1}^N c_k e_k  \nr_{p, \omega}^p =
\sum_{k=1}^N |c_k|^p .\ens Now by Proposition \ref{190609b} the stated results for the right-shift
operator $R$ follows immediately. It also follows that the left-shift operator
$L$ is bounded and that $||L|| \le 1.$  Since $||Le_2||_{p,w}=||e_1||_{p,w}=
||e_2||_{p,w},$ we finally conclude that $||L||=1,$ as claimed. \ep

While the scaling of the basis in $\{\delta_k\}_{k=1}^\infty$ in Corollary
\ref{190610a} indeed makes the operators $L$ and $R$ bounded, the scaling might
affect other conditions that are put on the basis, see., e.g., the condition (ii)
in the forthcoming Theorem \ref{050519a}. For the case that it is most convenient
to work with the canonical unit vector basis, we now characterize the weighted
$\ell^p$-spaces for which the left/right-shift operators are bounded
with respect to $\{\delta_k\}_{k=1}^\infty$. We leave the proof to the reader.

\bl\label{210519d} Fix any $ p\in[1,\infty)$, consider a sequence of positive scalars $\{w_k\}_{k=1}^\infty$, and let $\{\delta_k\}_{k=1}^\infty$
denote the standard basis for $\ell^p_w$. Then the following holds true:
\bei\item[(i)] The left-shift operator
$L\delta_1=0, \, L\delta_k=\delta_{k-1}. \, k\ge 2,$ extends  to a bounded
operator on $\ell^p_w$ if and only if
$\sup_{k\geq2} \frac{w_{k-1}}{w_k} <\infty; $ in the affirmative case, $$\| L \| = \sup_{k\geq 2} \left(\frac{w_{k-1}}{w_k} \right)^{1/p}.$$
\item[(ii)] The right-shift operator $R\delta_k=\delta_{k+1} $ extends
to a  bounded operator on $\ell^p_w$ if and only if
$\sup_{k\geq 2} \frac{w_{k}}{w_{k-1}}<\infty; $ in the affirmative case, $$\| R \| =\sup_{k\geq 2} \left(\frac{w_{k}}{w_{k-1}}\right)^{1/p}. $$
\eni\el

Several constructions of $p$-Riesz bases are available in the literature. Let us
mention an example in the setting of shift-invariant subspaces of $L^p(\mr)$:

\bex\rm It is well-known \cite{BL1} how to construct
a  function $\varphi \in \ltr$ such that the set of integer-translates
$\{\varphi (\cdot -k)\}_{k\in \mz}$ form a Riesz basis for the
Hilbert space
$ S_2:= \span \{\varphi (\cdot -k)\}_{k\in \mz}.$ Requiring furthermore
that $\varphi$
belongs to the
{\it Wiener space,} i.e., that $\sukz || \varphi \chi_{[k,k+1)
}||_\infty < \infty,$ it is proved in \cite{AST} that for
any $p\in [1, \infty)$ the family $\{\varphi (\cdot -k)\}_{k\in \mz}$
is a $p$-Riesz basis for the subspace $S_p$ of $L^p(\mr)$ given by
\bes S_p:= \left\{ \sukz c_k \varphi(\cdot -k)  \, \big|  \ \{c_k\}_{k\in \mz} \in \ell^p(\mz)\right\}.\ens Similar results are known
in the setting of modulation spaces,
introduced by Feichtinger in \cite{Fe4}.  \enx

Let us finally mention that the definition of $p$-Riesz bases can
be generalized in an obvious way to the so-called $X_d$-Riesz bases \cite{CT};
here, the sequence space $\ell^p$ is simply replaced by a general Banach sequence
space $X_d.$ Rather than going for the highest level of abstraction
we have decided to state the results in the setting of $p$-Riesz bases,
because this setting allows us to be very explicit and hereby facilitates
concrete applications.

\section{Approximate operator representations}
The goal of this section is to derive approximate representations
of sequences $\ftk$ in a separable Banach space $X$ of the
form $\{T^{\alpha(k)} \varphi\}_{k=1}^\infty$ for  suitable choices of
a vector $\varphi \in X$, a
bounded operator $T: X \to X,$ and the powers $\alpha(k), \, k\in \mn.$
We begin with the case where $\ftk$ consists of ``finite sequences"
in Section \ref{181205a}. The results are generalized to ``sufficiently fast decaying
sequences" in Section \ref{181205b}. An alternative approach, tailored to
the setting of Banach function spaces, is presented in Section \ref{190603a}.
The purpose of Sections \ref{181205c}--\ref{181205d}
is to apply the results in Sections \ref{181205a}--\ref{190603a} to
the setting of atomic decompositions and Banach frames: we
show how to design the approximations such that the sequence
$\{T^{\alpha(k)} \varphi\}_{k=1}^\infty$ keeps essential features of the
given sequence $\ftk.$

In order to facilitate reading of the next sections, we mention that the theoretical
results in Sections \ref{181205a}--\ref{190603a} have a common structure: fixing
arbitrary positive scalars $\{\epsilon_k\}_{k=1}^\infty\in\ell^1$ that are chosen according to
the desired level of approximation, they show how we for a given sequence
$\ftk \subset X$ can choose a vector $\varphi \in X$, a
bounded operator $T: X \to X,$  and corresponding  powers $\alpha(k), \, k\in \mn,$
such that
\bee\label{111218cn} \| f_k - T^{\alpha(k)}\varphi \| \le \sum_{j=k+1}^\infty \epsilon_j. \ene Since the positive scalars $\{\epsilon_k\}_{k=1}^\infty\in\ell^1$ are arbitrary,
this implies that we for any given sequence $\{{\mathcal E}_k\}_{k=1}^\infty$ of positive
scalars can obtain that $\| f_k - T^{\alpha(k)}\varphi \| \le {\mathcal E}_k.$
Indeed, assuming without loss of generality that the sequence
$\{{\mathcal E}_k\}_{k=1}^\infty$ is decreasing and that ${\mathcal E}_k\to 0$
as $k\to \infty,$ we obtain this inequality by applying
\eqref{111218cn} to any sequence $\{\epsilon_k\}_{k=1}^\infty\in\ell^1$ such that
$\epsilon_k \le {\mathcal E}_k-{\mathcal E}_{k+1}$ for all $k\in \mn.$
The consequences of an equality of the
form  \eqref{111218cn} are considered in  Sections \ref{181205c}--\ref{181205d}.

\subsection{Finite sequences in Banach spaces} \label{181205a}

Let $X$ denote a separable Banach space having a basis $\etk$. As standing assumption we need that the operators $L,R$ defined in \eqref{050519b} extend to bounded linear operators on $X$. Choose $\lambda>\| R\| $ and consider the bounded operators
$T, S: X \to X$ given by
\bee\label{050519d}
T=\lambda L,\quad S=\lambda^{-1} R
\ene
Note that by the choice of $\lambda$ we have $\| S\| <1$, a condition that will be crucial in Theorem \ref{18314b}.  We will first consider natural generalizations of finite sequences
to the general Banach space $X.$ Indeed, we will assume that each vector $f_k$ is a
finite linear combination of vectors from the basis $\etk$. Our approach and
proof are inspired by a classical result by Rolewicz \cite{rolewicz} in the
setting of hypercyclic operators.

\bt \label{18314b} In the above setup, consider a sequence $\ftk\subset X$ and assume that for every $k\in\mn$, there
exists an integer $N(k)$ such that $f_k \in \Span \{e_j\}_{j=1}^{N(k)}$.
Fixing any sequence $\{\epsilon_k\}_{k=1}^\infty\in\ell^1$ consisting of positive numbers,
choose a sequence $\{\alpha(k)\}_{k\in\mn}$ of nonnegative integers   such that
\bee\label{111218b} \qquad\alpha(k)- \alpha(k-1) \ge \max\left\{~\frac{ \ln\epsilon_k-\ln{\|f_k\|}}{\ln\|S\|},~ N(k-1), N(k)~\right\}, \quad k\ge 2,\ene
and let \bee \label{19314a} \varphi:=\sum_{j=1}^{\infty} S^{\alpha(j)} f_j.\ene Then
\bee\label{111218c} \| f_k - T^{\alpha(k)}\varphi \| \le \sum_{j=k+1}^\infty \epsilon_j. \ene
\et
\bp For convenience of the proof, let
\bee\label{111218a}
r_k := \max\left\{~\frac{ \ln\epsilon_k-\ln{\|f_k\|}}{\ln\|S\|},~ N(k)~\right\}, \quad k\in\mn.
\ene
We first show that the vector $\varphi$ in \eqref{19314a} is well-defined. To that end, let $\varphi_n:=\sum_{j=1}^{n} S^{\alpha(j)} f_j, \, n\in \mn.$
Then,  for any $m,n\in\mn$ with $m\leq n$, we have
\bes \| \varphi_n - \varphi_m\|=
\| \sum_{j=m}^{n} S^{\alpha(j)} f_j \| \leq \sum_{j=m}^{n} \|S\|^{\alpha(j)}\| f_j \|.\ens
It follows from \eqref{111218a} that
$\|S\|^{r_j}\|f_j\| \leq \epsilon_j$ for all $j\in \mn.$ Since $\| S \| <1 $ and
$\alpha(k)>r_k$ by \eqref{111218b} and the definition of $r_k,$ it follows that
\bes \| \varphi_n - \varphi_m\| \leq \sum_{j=m}^{n} \|S\|^{r_j}\| f_j \|
\leq \sum_{j=m}^{n} \epsilon_j \to 0\quad\mbox{ as }m,n\to\infty. \ens
Thus $\{\varphi_n\}_{n=1}^\infty$ is a Cauchy sequence  and hence $\varphi$
is well-defined. We now prove that \eqref{111218c} holds.
Fix $k\in\mn$ and consider $j\in\{1,\dots,k-1\}$.  The inequality \eqref{111218b} implies that \bes \alpha(k)-\alpha(j)> \alpha(j+1)-\alpha(j) \ge N(j). \ens
Therefore for any $j< k$, we have $ T^{\alpha(k)-\alpha(j)}f_j  = 0. $
Thus we can write
\bes \| f_k - T^{\alpha(k)}\varphi \| &=&  \| f_k - T^{\alpha(k)}\sum_{j=1}^{\infty} S^{\alpha(j)} f_j \|\\
&=&\| \sum_{j=1}^{k-1}T^{\alpha(k)-\alpha(j)}f_j +
\sum_{j=k+1}^{\infty}S^{\alpha(j)-\alpha(k)}f_j \|\\
&=& \|\sum_{j=k+1}^{\infty}S^{\alpha(j)-\alpha(k)}f_j \|\\
&\leq& \sum_{j=k+1}^{\infty}\|S\|^{\alpha(j)-\alpha(k)} \|f_j \|.\ens
Now for $j>k$,  \eqref{111218b} implies
$$ \alpha(j)-\alpha(k)\ge \alpha(j)-\alpha(j-1)\ge r_j. $$ Therefore
\bes
\| f_k - T^{\alpha(k)}\varphi \| \leq \sum_{j=k+1}^{\infty}\|S\|^{r_j}\|f_j \| \leq \sum_{j=k+1}^{\infty}\epsilon_j,
\ens  as claimed.\ep

%

The key condition in Theorem \ref{18314b}   is that the left/right-shift operators
with respect to a certain basis are bounded. Recall that
Corollary \ref{190610a} shows how to fulfill this condition in any
weighted $\ell^p$-space, $1\le p < \infty.$
It is typically significantly more complicated to verify boundedness of
the shift operators on Banach function spaces than on Banach sequence spaces.
For this reason we will formulate an alternative result in Section
\ref{190603a}, tailored to the setting of Banach
function spaces.

\subsection{Localized sequences in Banach spaces} \label{181205b}
Next we will  prove that Theorem \ref{18314b} can be generalized to certain
infinite sequences, provided their coordinates decay ``sufficiently fast".
To motivate the exact formulation of the condition, assume for a moment that
$\ftk$ is a sequence in $\ell^p$ for some $p\in (1, \infty).$ The canonical delta-sequence $\etk$ is an unconditional basis for
$\ell^p$ and its dual space $\left(\ell^p\right)^*$, and the $jth$ coordinate
in the vector $f_k$ is precisely $\la f_k, e_j\ra.$ A natural way of defining
``fast decay" of the coordinates of $f_k$ is to require that there exist
constants $C,\beta>0$ such that
\bee \label{041018a} | \la f_k,e_j\ra| \le C e^{-\beta \, |j-k|}, \,\quad \forall  j,k \in \mn.\ene
We will use exactly this idea, but formulated for a general basis for the Banach space $X$.

\bt\label{050519a} Let $X$ denote a  Banach space with basis $\etk$
and associated dual basis $\{e_k^*\}_{k=1}^\infty$, and let $X_d$ be a solid Banach sequence space with an absolutely continuous norm, which
contains the canonical unit vectors $\{\delta_k\}_{k=1}^\infty.$ Let $\ftk\subset X$ and assume the followings:  \bei
\item[(i)] The left/right-shift operators $L,R$ with respect to the given basis $\etk$ are bounded on $X$. Choose any $\lambda> ||R||$.
\item[(ii)] There exists a constant $B>0$ such that
\[  \| \suk c_k e_k \| \leq B \| \ctk \|_{X_d},  \]
for all finite sequences $\ctk$.
\item[(iii)] There exist constants $C>0$ and $\beta > \ln \lambda $ such that
$\{ e^{-\beta j}\}_{j=1}^\infty\in X_d$ and
\bee\label{130619a} | \la f_k , e_j^* \ra | \leq C e^{-\beta |j-k|} \quad\forall  j,k\in\mn. \ene
\eni
Finally, fixing a sequence $\{\epsilon_k\}_{k=1}^\infty\in\ell^1$ of positive
scalars, choose a strictly increasing sequence of nonnegative integers $\{\alpha(k)\}_{k=1}^\infty$  such that for all $k\in\mn$,
\bee\label{180319aa}
\| S \|^{\alpha(k)-\alpha(k-1)} \| f_k \|_X \leq \epsilon_k/2.\ene
Then $\varphi := \sum_{k=1}^\infty S^{\alpha(k)}f_k$ is well-defined.
Moreover, by choosing $\alpha(1)=0$ and $\{\alpha(k)\}_{k=2}^\infty$ recursively such that
\bee  \label{180319cc}
\alpha(k)  \geq  \alpha(k-1) +  k - 2,\quad \forall  k\geq 2,
\ene
and
\bee \notag
\alpha(k)
& >  & \frac{\ln( \sum\limits_{j=k+1}^{\infty}  \epsilon_j)-\ln(\left\|\{ e^{-\beta j} \}_{j=1}^\infty \right\|_{X_d})-\ln(\sum\limits_{n=0}^{k-1}  (\lambda e^{-\beta})^{-\alpha(n)}  e^{\beta n}) -
	\ln(2BC)}{\ln(\lambda) -\beta}, \\
& \, &  \label{200319cc} \ene
then
\bes \| f_k - T^{\alpha(k)}\varphi \|_X \le \sum_{n=k+1}^{\infty}  \epsilon_n, \,\quad \forall  k\in \mn.\ens

\et
\bp First note that the infinite sum $\sum_{k=1}^\infty S^{\alpha(k)}f_k$ is absolutely convergent; indeed, by \eqref{180319aa},
\[ \sum_{k=1}^\infty \| S^{\alpha(k)}f_k \|_X \leq \suk \| S \|^{\alpha(k)} \| f_k \|_X
\leq  1/2\sum_{k=1}^\infty \epsilon_k <\infty.   \]
Thus $\varphi = \sum_k S^{\alpha(k)}f_k$ is well-defined.
Let $k\in\mn$, then
\bee \label{192706a}
\qquad\| f_k - T^{\alpha(k)}\varphi \|_X \leq   \sum_{n=1}^{k-1} \|T^{\alpha(k)}S^{\alpha(n)}f_n \|_X + \| \sum_{n=k+1}^{\infty} T^{\alpha(k)}S^{\alpha(n)}f_n \|_X.
\ene
We study the two terms at the right-hand side of the inequality separately.
First,
\bes
\| \sum_{n=k+1}^{\infty} T^{\alpha(k)}S^{\alpha(n)}f_n \|_X & = &
\| \sum_{n=k+1}^{\infty} S^{\alpha(n)-\alpha(k)} f_n \|_X \leq     \sum_{n=k+1}^{\infty} \|S^{\alpha(n)-\alpha(k)} f_n \|_X\nonumber\\
&\leq &    \sum_{n=k+1}^{\infty} \|S\|^{\alpha(n)-\alpha(k)} \|f_n \|_X\nonumber\\
&\leq& \sum_{n=k+1}^{\infty} \| S\|^{\alpha(n)-\alpha(n-1)} \| f_n\|_X.\ens
Using \eqref{180319aa}, we get
\bee\label{180319dd} \| \sum_{n=k+1}^{\infty} T^{\alpha(k)}S^{\alpha(n)}f_n \|_X  \leq 1/2\sum_{n=k+1}^{\infty}  \epsilon_n.
\ene
Now, for $n=1,\dots,k-1$,

\bes
\|T^{\alpha(k)}S^{\alpha(n)}f_n \|_X &=& \| T^{\alpha(k)} S^{\alpha(n)} \sum_{j=1}^\infty \la f_n , e_j^* \ra e_j \| \\
&=& \|  \lambda^{\alpha(k)-\alpha(n)}  \sum_{j=\alpha(k)-\alpha(n)+1}^\infty \la f_n , e_j^* \ra e_{j-\alpha(k)+\alpha(n)}             \| \\
&=& \|  \lambda^{\alpha(k)-\alpha(n)}  \sum_{j=1}^\infty \la f_n , e_{j+\alpha(k)-\alpha(n)}^* \ra e_{j}             \|. \\
\ens
Using the condition (iii), we have $| \la f_n , e_{j+\alpha(k)-\alpha(n)}^*\ra|\leq Ce^{-\beta|j+\alpha(k)-\alpha(n)-n|}$;
since $X_d$ is a solid Banach sequence space it follows that $\{\la f_n , e_{j+\alpha(k)-\alpha(n)}^* \ra\}_{k=1}^\infty \in X_d$ and that
\bes  \left\| \{ \la f_n , e_{j+\alpha(k)-\alpha(n)}^*\ra \}_{j=1}^\infty \right\|_{X_d}\leq C\left\|\{ e^{-\beta |j+\alpha(k)-\alpha(n)-n	| } \}_{j=1}^\infty \right\|_{X_d}.\ens
A standard argument (using
that $X_d$ is assumed to have an absolutely continuous norm) shows that
the condition (ii) actually implies that the stated inequality holds
for all $\ctk \in X_d;$ thus we arrive at
\bes\label{180319ee}
\| T^{\alpha(k)}S^{\alpha(n)}f_n \|_X &\leq&  BC \lambda^{\alpha(k)-\alpha(n)} \left\|\{ e^{-\beta |j+\alpha(k)-\alpha(n)-n	| } \}_{j=1}^\infty \right\|_{X_d}.   \ens
Condition \eqref{180319cc} implies that $j+\alpha(k)-\alpha(n)-n\geq 0$. Thus
\bes
\| T^{\alpha(k)}S^{\alpha(n)}f_n \|_X& \leq  &  BC \lambda^{\alpha(k)-\alpha(n)}
\left\|\{ e^{-\beta (j+\alpha(k)-\alpha(n)-n)	} \}_{j=1}^\infty \right\|_{X_d}\\
&=& BC \lambda^{\alpha(k)-\alpha(n)}
\left\| e^{-\beta (\alpha(k)-\alpha(n)-n)} \{e^{-\beta j	} \}_{j=1}^\infty \right\|_{X_d}\\
&\leq&
BC (\lambda e^{-\beta})^{\alpha(k)-\alpha(n)}  e^{\beta n}
\left\|\{ e^{-\beta j} \}_{j=1}^\infty \right\|_{X_d}.
\ens
If $\{\alpha(k)\}_{k=1}^\infty$ satisfies the growth condition specified in \eqref{200319cc}, then we conclude
\bee\label{290419a}
\sum_{n=1}^{k-1}\| T^{\alpha(k)}S^{\alpha(n)}f_n \|_X \leq 1/2\sum_{j=k+1}^\infty \epsilon_j.
\ene
The result now follows from \eqref{192706a},  \eqref{180319dd}   and \eqref{180319ee}. \ep

\bex \rm Consider the Banach space $X=\ell^p_w$. As we proved in  Corollary \ref{190610a}, the right/left-shift operators $L,R$ defined with respect to the basis $\etk:= \{ w_k^{-1/p}\delta_k\}_{k=1}^\infty$ are bounded and $\|L\|=\|R\|=1$. Let $X_d=\ell^p$. Clearly $\ell^p$ is a solid Banach space with absolutely continuous norm and it contains the canonical basis $\dtk$. Moreover, for every finite sequence $\{c_k\}$,
\[ \|\sum c_k e_k \|_X = \| \{c_k w_k^{-1/p}\} \|_X =  \left(\sum |c_k|^p\right)^{1/p}=\| \{c_k\} \|_{X_d}.\]
Therefore all the conditions  in Theorem \ref{050519a} are satisfied. The dual basis of $\etk$ is given by $e_k^*=w_k^{1/p}\delta_k$, $k\in\mn$. Given any sequence $\ftk\subset X$, as in Theorem \ref{050519a}, write $f_k=\{(f_k)_j \}_{j=1}^\infty$. Then
\[  |\la f_k , e_j^* \ra | = |  (f_k)_j | w_j^{1/p}, \, \,  j,k\in\mn. \]
This shows that in the setting of $X=\ell^p_w, X_d=\ell^p$, Theorem \ref{050519a} applies to all sequences $\ftk$ such that for some $C,\beta>0$,
\bes |  (f_k)_j |   \leq Cw_j^{-1/p}e^{-\beta|j-k|},\,\quad\forall  j,k\in\mn. \ens 
\enx

\subsection{Banach function spaces} \label{190603a}
The results in Sections \ref{181205a}--\ref{181205b} deal with general Banach spaces, having a basis with respect to which the shift operators $L$ and $R$ are bounded.
This condition is often considerably more complicated to verify in
Banach function spaces than in Banach sequence spaces.
For this reason we will now consider the special case of Banach function spaces, but without
any condition of knowledge of a basis such that the corresponding left/right-shift operators are bounded.
In the entire section we let  $X$ denote a Banach space of functions $f:\mr\to\mc$.
For $a\in \mr,$ consider the translation operator $T_a$ acting on functions
$f:\mr \to \mc$ by $T_af(x):=f(x-a).$ We say that $X$ is {\it translation invariant}
if the translation operators $T_1$ and $T_{-1}$ map functions in $X$ into $X;$
in this case, if $T_1$ is bounded, it follows from the open mapping theorem
that also $T_{-1}$ is bounded.
Assuming that $X$ is a solid Banach space and the translation operator $T_1$
is bounded, consider now for any $\lambda>\| T_1\|$   the weighted translation operators
$S, T: X \to X,$
\bee \label{190610c}  Tf=\lambda (T_{-1} f)\chi_{[0,\infty)},\quad  Sf=\lambda^{-1}T_1f,\quad  f\in X.\ene
Note that $\| S\|<1 $ and
that for $k\in \mn,$ we have $T^k f= \lambda^k (T_{-k}f) \chi_{[0,\infty)}, \, f\in X.$

Before stating the main result, Theorem \ref{190610d}, let us comment on one of
the conditions in the statement of the result and various ways
of circumventing it. We will consider a sequence $\ftk$ of function in $X$ that are supported in $[0,\infty)$. First, the choice of the interval $[0,\infty)$ is
not essential: the result immediately generalizes to functions supported
on any  half interval $[a,\infty), a\in\mr$, simply by replacing
the characteristic function $\chi_{[0,\infty)}$ in  the translation operator $T$ defined in \eqref{190610c}  by $\chi_{[a,\infty)}$. Next, if the sequence $\ftk$ can be written as $\ftk=\gtk \cup \htk$ where $ \supp\, g_k\subset [0,\infty)$ and $\supp\,h_k\in (-\infty,L]$ for some $L\in\mr$, then a similar procedure as suggested in the
particular case of a Hilbert space in \cite{olemarzieh-10} can be applied on $\gtk$ and $\htk$ separately. In this case, the sequence $\gtk \cup \htk$ can be approximated with a union of two suborbits, each
associated with a bounded operator. We refer the interested reader to \cite{olemarzieh-10}
for details.

\bt \label{190610d} Let $X$ denote a solid translation-invariant Banach
function space with absolutely
continuous norm, and assume that the translation operator $T_1$  acts
boundedly on $X.$
Let $\lambda >\| T_1\|$ and $\mu> \lambda\|T_{-1}\|$. Assume that $\ftk\subset X$ and $\mathrm{supp} f_k\subset [0,\infty)$. Also assume that for every $k\in\mn$, there exist $a_k\in\mn$  and $C_k>0$ such that for every $a\geq a_k$,
\bee \label{240519a}
\| f_k\chi_{[a,\infty)} \| \leq C_k\mu^{-a}.
\ene
Fixing now a sequence $\{\epsilon_k\}_{k=1}^\infty\in\ell^1$
of positive scalars, choose an increasing sequence $\{\alpha(k)\}_{k=1}^\infty$ of nonnegative integers such that
\bee\label{240519b}
\| S \| ^{\alpha(k)-\alpha(k-1)}\| f_k \| < 1/2\epsilon_k,\quad  k\geq 2.
\ene
Then $\varphi:=\suk S^{\alpha(k)}f_k$ is well-defined.
Furthermore, if we also assume that $\alpha(1)=0$ and
\bee\label{240519c}
\alpha(k+1)-\alpha(k)\geq a_{k},\quad  k\in\mn,
\ene
and
\bee\label{240519d}
\qquad\quad\alpha(k) \geq \frac{\ln(2)+ \ln\left( \sum_{n=1}^{k-1}\,C_n\Big((\lambda \|T_{-1}\|)^{-1}\mu\Big)^{\alpha(n)}\right)-\ln(\sum_{j=k+1}^\infty \epsilon_j)}{\ln((\lambda \| T_{-1}\|)^{-1}\mu)},
\ene
then \bes \| f_k - T^{\alpha(k)}\varphi \| \leq \sum_{n=k+1}^\infty \epsilon_n.\ens
\et
\bp For every $m,n\in\mn$, by \eqref{240519b} and since $\|S\|<1$,
\bes \| \sum_{k=m}^N S^{\alpha(k)} f_k \| &\leq&  \sum_{k=m}^N \| S\|^{\alpha(k)}\| f_k \| \\
&<& 1/2\sum_{k=m}^N \epsilon_k\to 0\quad\text{ as }m,n\to\infty.
\ens
As in the proof of Theorem \ref{18314b}, it follows that $\varphi$ is well-defined. Also for $k\in\mn$, we have
\bee \label{192706b}
\| f_k - T^{\alpha(k)}\varphi \| \leq   \sum_{n=1}^{k-1} \|T^{\alpha(k)}S^{\alpha(n)}f_n \| + \| \sum_{n=k+1}^{\infty} T^{\alpha(k)}S^{\alpha(n)}f_n \|.
\ene
We now consider the two terms at the right-hand side of the inequality separately.
First, using \eqref{240519b}, we obtain
\bes
\| \sum_{n=k+1}^{\infty} T^{\alpha(k)}S^{\alpha(n)}f_n \| & = &
\| \sum_{n=k+1}^{\infty} S^{\alpha(n)-\alpha(k)} f_n \| \leq     \sum_{n=k+1}^{\infty} \|S^{\alpha(n)-\alpha(k)} f_n \|\nonumber\\
&\leq &    \sum_{n=k+1}^{\infty} \|S\|^{\alpha(n)-\alpha(k)} \|f_n \|
\leq \sum_{n=k+1}^{\infty} \| S\|^{\alpha(n)-\alpha(n-1)} \| f_n\|\\
&\leq & 1/2\sum_{n=k+1}^{\infty}  \epsilon_n.
\ens
Next, we get

\bes
\sum_{n=1}^{k-1}\|T^{\alpha(k)}S^{\alpha(n)}f_n \|&=& \sum_{n=1}^{k-1}\|T^{\alpha(k)-\alpha(n)}f_n \| \nonumber \\
&=&\sum_{n=1}^{k-1}\|\lambda ^{\alpha(k)-\alpha(n)}(T_{-1}^{\alpha(k)-\alpha(n)}f_n)\chi_{[0,\infty)} \| \nonumber\\
&=& \sum_{n=1}^{k-1}\lambda^{\alpha(k)-\alpha(n)}\| T_{-1}^{\alpha(k)-\alpha(n)}\big(f_n\chi_{[\alpha(k)-\alpha(n),\infty)} \big)\| \nonumber\\
&\leq& \sum_{n=1}^{k-1} (\lambda \| T_{-1}\|)^{\alpha(k)-\alpha(n)}\,\| f_n\chi_{[\alpha(k)-\alpha(n),\infty)} \|
.\ens
Since $n\leq k-1$, we have $\alpha(k)-\alpha(n) \geq \alpha(n+1)-\alpha(n)$ and therefore by \eqref{240519c}, we get
$\alpha(k)-\alpha(n) \geq a_n$. Therefore, by \eqref{240519a},
\bee\label{260519a}
\sum_{n=1}^{k-1}\|T^{\alpha(k)}S^{\alpha(n)}f_n \|&\leq &  \sum_{n=1}^{k-1}(\lambda \| T_{-1}\|)^{\alpha(k)-\alpha(n)} C_n\mu^{-(\alpha(k)-\alpha(n))}\nonumber\\
&=& (\lambda \| T_{-1}\| \mu^{-1})^{\alpha(k)}  \sum_{n=1}^{k-1} C_n((\lambda \| T_{-1}\|)^{-1} \mu)^{\alpha(n)}.
\ene
Now if $\alpha(k)$ is chosen such that \eqref{240519d} holds,  we have that
\bes
\alpha(k)\ln\Big((\lambda \| T_{-1}\|)^{-1}\mu\Big) \geq \ln(2)+ \ln\left( \sum_{n=1}^{k-1}\,C_n((\lambda \|T_{-1}\|)^{-1}\mu)^{\alpha(n)}\right)-\ln\Big(\sum_{j=k+1}^\infty \epsilon_j\Big)
\ens
or
\bes
\alpha(k)\ln\Big(\lambda \| T_{-1}\| \mu^{-1}\Big)+  \ln\left( \sum_{n=1}^{k-1}\,C_n((\lambda \|T_{-1}\|)^{-1}\mu)^{\alpha(n)}\right) \leq\ln\Big(1/2\sum_{j=k+1}^\infty \epsilon_j\Big).
\ens
Applying the exponential function on both sides of the inequality yields
\bes
(\lambda \| T_{-1}\| \mu^{-1})^{\alpha(k)}  \sum_{n=1}^{k-1}\, C_n( (\lambda \|T_{-1}\|)^{-1}\mu)^{\alpha(n)} \leq 1/2\sum_{j=k+1}^\infty \epsilon_j,
\ens
and from \eqref{260519a}, we conclude that
$$\sum_{n=1}^{k-1}\|T^{\alpha(k)}S^{\alpha(n)}f_n \|\leq  1/2 \sum_{j=k+1}^{\infty} \epsilon_j.$$ Using now \eqref{192706b} and the obtained estimates of the
two terms, we obtain the desired result.
\ep

In the next example we consider an important class of Banach spaces
that satisfy all the conditions in Theorem \ref{190610d}.

\bex\rm
Let $m:\mr\to[0,\infty)$ be a continuous function and  $w:\mr\to[0,\infty)$  a $m$-moderate weight function, i.e., a measurable function such that
$$ w(x+y)\leq m(x)w(y), \quad\forall x,y\in\mr.$$
For $1\leq p<\infty,$ let
\bes L^p_w(\mr):= \left\{f: \mr \to \mc \, \big| \,\int_{\mr} |f(x)|^p w(x) dx < \infty \right\}.  \ens
Then $L^p_w(\mr)$ is a Banach space with respect to the norm
\[ \| f\|_{L^p_w}=\left( \int_{\mr} |f(x)|^p w(x) dx \right)^{1/p}.
\]
We leave it to the reader to verify that the norm is absolutely continuous, $L^p_w(\mr)$ is invariant under the translation-operators $T_1,T_{-1},$ and that
$||T_{-1}|| \le m(-1)^{1/p}  $ and $||T_1|| \le m(1)^{1/p}$.
 \enx

In order to demonstrate the practical issues showing up in applications of
Theorem \ref{190610d}, we will now consider so-called {\it Gabor systems.}
For $b\in\mr$, let $E_{b}:\ltr\to\ltr$ be the modulation operator defined as $E_b g(x)=e^{2\pi i bx}g(x)$, $x\in\mr$. Fixing a function
$g\in \ltr$ and some parameters $a,b>0,$  the sequence of functions $\{E_{mb}T_{na} g\}_{m,n\in\mz}$ in $\ltr$ is called a Gabor system. Since the translation
operators and modulation operators clearly act boundedly on any $L^p$-space, $1\le p < \infty$ as well, we will now assume that $g\in L^p(\mr)$ and consider the
Gabor system in $L^p(\mr)$ instead. Note that if the function $g$ is compactly
supported, we can split the Gabor system into a union $\gtk \cup \htk$ where $ \supp\, g_k\subset [0,\infty)$ and $\supp\,h_k\in (-\infty,L]$ for some $L\in\mr;$ thus,
as explained just before the statement of Theorem \ref{190610d} we can
approximate the Gabor system using suborbits of two bounded operators.
In the next example we will replace the assumption that $g$ has compact support by
the assumption that $\supp\, g\subset [0,\infty),$
and show how to obtain the estimate \eqref{240519a} for a certain ordering of a
the ``half Gabor system"  $\{E_{mb}T_{na} g\}_{m\in\mz,n\in\mn\cup\{0\}}$.

\bex\rm Consider a function $g\in L^p(\mr), \, 1\le p < \infty,$ and
assume that $\supp\, g\subset [0,\infty)$ and  that there exist constants $C,d_0>0$ and $\mu>1$ such that for all $d \geq d_0$
\bee \label{030619a}
\left(	\int_{d}^\infty | g(x) |^p dx \right)^{1/p} \leq C\mu^{-d}.
\ene
Re-index the ``half Gabor system" $\{E_{mb}T_{na} g\}_{m\in\mz,n\in\mn\cup\{0\}}$ as $\ftk$  in such a way that  $f_k$ and $f_{k+1}$ differ with  at most one translate by $a$, i.e., if $f_k = E_{mb}T_{na} g$ for some $m\in\mz$ and $n\in\mn$ then $f_{k+1}$ is one of the following functions
\[  E_{mb}T_{(n\pm1)a} g, \quad E_{(m\pm 1)b} T_{na}g.\]
Now, for $k\in\mn$, write $f_k= E_{mb}T_{na} g,$ where $n\in\{0,1,\dots,k\}$ and $m\in\mz$. Then, considering any $d>0$, we have
\bes
\int_{d}^{\infty} | f_k(x)|^p dx &=&
\int_{d}^{\infty} | E_{mb}T_{na}g(x)|^p dx
= \int_{d}^{\infty} | g(x-na)|^p dx = \int_{d-na}^{\infty} | g(x)|^p dx. \ens
For $d\geq d_0 + ka,$  using that $n\in \{0,1,\dots,k\}$, \eqref{030619a} yields
\bes
\left(\int_{d}^{\infty} | f_k(x)|^p dx \right)^{1/p}&\leq& C\mu^{-(d'-na)}
\leq  (C\mu^{ka}) \mu^{-d}.
\ens
Thus, choosing $C_k:= (C\mu^{ka}) $ and $d_k:=d_0+ka$, we have that  for any $d\geq d_k$,
\bee\label{030619b} \| f_k \chi_{[d,\infty)}\|_p \leq C_k\mu^{-d}, \ene
i.e., the condition \eqref{030619a} is satisfied.  By a direct calculation, the ``new technical assumption"
\eqref{030619a} holds if, e.g.,
\bes |g(x)| \leq e^{-x} \chi_{[0,\infty)}, \quad\forall x\in\mr;\ens
indeed, in this  we can take
$C=1/p$ and $\mu=e^p$. Note that for $p=2$, the function $h(x)=e^{-x} \chi_{[0,\infty)} $
play a special role in Gabor analysis: it generates a Gabor frame $\{E_{mb}T_{na}h\}_{m,n\in\mz}$ for $L^2(\mr)$ if and only if $ab\leq 1$, see \cite{Janssen}. 
\enx

\subsection{ $\epsilon$-close approximations of $\ftk$     } \label{181205c}

In this section we will pave the way for the application of
the theoretical results in Section \ref{181205d}.
To motivate what follows, consider again
Theorem \ref{18314b}:  for any given finite sequence
$\ftk$ in $X$ and any sequence $\{\epsilon_k\}_{k=1}^\infty\in\ell^1$ of positive scalars it
specifies how to choose $\varphi\in X$ and powers $\alpha(k)$ such that
for each $k\in \mn$ the vector
$T^{\alpha(k)}\varphi$ belongs to a ball  around $f_k,$ with a radius specified
by \eqref{111218c}.
The goal of typical approximation arguments is that the stated condition should
imply that the approximating sequence - here $\{ T^{\alpha(k)}\varphi\}_{k=1}^\infty$ -
shares key features of the given sequence $\ftk.$ In concrete cases this might call for extra conditions
on the sequence $\{\epsilon_k\}_{k=1}^\infty\in\ell^1$. Recall (see, e.g., \cite{Y}) that given a sequence
$\ftk$ in a Hilbert space $\h,$ a sequence $\gtk \subset \h$ is said to be
{\it quadratically close} to $\ftk$ if
\bes \sum_{k=1}^\infty ||f_k-g_k||^2 < 1.\ens
This concept is well-motivated. Indeed, if $\ftk$ is an orthonormal basis for $\h$
and $\gtk \subset \h$ is quadratically close to $\ftk,$ then $\gtk$ is a Riesz basis
for $\h;$ see \cite{Y}. The following generalizes  the concept of quadratically
close sequences.




\bd Let $\ftk $ be a sequence in a Banach space $X.$
	Given any $p\in [1, \infty)$ and any $\epsilon >0,$ a sequence $\gtk\subset X$ is said to be $\epsilon$-close to
	$\ftk$ with respect to $\ell^p$ if
	\bee \sum_{k=1}^\infty ||f_k-g_k||^p \le \epsilon.\ene \ed

In the setting of Theorem \ref{18314b}, Theorem \ref{050519a} or
Theorem \ref{190610d} we will now
show how to ensure that the constructed sequence $\{ T^{\alpha(k)}\varphi\}_{k=1}^\infty$ is $\epsilon$-close to the given sequence $\ftk$. Recall that
on the structural level, the three mentioned results are similar: all
of them ensure that for a given sequence $\ftk$
in the considered Banach space and a fixed sequence$\{\epsilon_j\}_{j=1}^\infty \in \ell^1,$
the operator $T,$ the vector $\varphi,$ and the
powers $\alpha(k)$ satisfy the inequality \eqref{111218cn}.

\bc  \label{18315a} Let $\ftk\subset X.$ Fix a positive sequence $\{ \epsilon_j \}_{j=1}^\infty\in\ell^1$, and assume that  the operator $T:X \to X,\, \varphi \in X,$
and powers $\alpha(k), \, k\in \mn,$ have been constructed such that
\eqref{111218cn} is satisfied.
Then the following holds true:

\bei
\item[(i)]  Let $\epsilon_j:= \epsilon 2^{-j}$. Then the sequence
$\{ T^{\alpha(k)}\varphi\}_{k=1}^\infty$ is $\epsilon$-close to $\ftk$ simultaneously
for all $p\in [1, \infty)$. More precisely, it holds that for all $p\geq 1$,
\bes \sum_{k=1}^\infty \| f_k - T^{\alpha(k)}\varphi \|^p \le \frac{\epsilon^p}{2^p-1} \le \epsilon.\ens
\item[(ii)] Fix any $p\in [1, \infty)$ and consider a weight sequence $\{w_k\}_{k=1}^\infty$ such that $\sum_{k=1}^\infty
\frac{1}{2^{kp}}  w_k<\infty$.  Let $M:=\max\{1,\sum_{k=1}^\infty
\frac{1}{2^{kp}}  w_k\}$ and  $\epsilon_j:=\frac{\epsilon}{M}2^{-j}$. Then
\bes \sum_{k=1}^\infty \| f_k - T^{\alpha(k)}\varphi \|^p w_k\leq \epsilon^p.\ens
\item[(iii)] Consider a solid Banach sequence space $X_d$
for which $\{2^{-k}\}_{k=1}^\infty\in X_d$. Take $M:=\| \{2^{-k}\}_{k=1}^\infty \|_{X_d}$ and let  $\epsilon_j:=\epsilon/M 2^{-j}$. Then
\bee\label{150319a} || \{ \| f_k - T^{\alpha(k)}\varphi\| \}_{k=1}^\infty ||_{X_d} \le \epsilon.\ene

\eni
\ec

\bp For the proof of (i), just observe that

\bes \sum_{k=1}^\infty \| f_k - T^{\alpha(k)}\varphi \|^p &\le&\sum_{k=1}^\infty \left( \sum_{j=k+1}^\infty \frac{ \epsilon}{2^j} \right)^p
=\epsilon^p \sum_{k=1}^\infty   \left( \frac{1}{2^{k}} \right)^p  = \frac{\epsilon^p}{2^p-1}.\ens
Similarly, under the assumptions in (ii),
\bes \sum_{k=1}^\infty \| f_k - T^{\alpha(k)}\varphi \|^p w_k &\leq &
\sum_{k=1}^\infty \left( \sum_{j=k+1}^\infty \epsilon/M 2^{-j} \right)^p w_k\\
& = &\sum_{k=1}^\infty  (\epsilon/M)^p
\frac{1}{2^{kp}}  w_k \\\\
&\leq&  \epsilon^pM^ {1-p} \leq \epsilon^p.\ens
For the proof of (iii), Theorem \ref{18314b} gives
\bee\label{200319a}
\| f_k - T^{\alpha(k)}\varphi \| \leq \epsilon /M2^{-k}.
\ene
Since $X_d$ is a solid Banach sequence space and $\{ 2^{-k}\}_{k=1}^\infty\in X_d$, this implies  $\{ \| f_k - T^{\alpha(k)}\varphi\| \}_{k=1}^\infty\in X_d$ and that
\bes
|| \{ \| f_k - T^{\alpha(k)}\varphi\| \}_{k=1}^\infty ||_{X_d} \leq \| \{\epsilon /M2^{-k} \}_{k=1}^\infty \|_{X_d} =\epsilon,
\ens as claimed.
\ep

\subsection{Applications to atomic decompositions}
\label{181205d}

We will now apply the theoretical results to atomic decompositions in Banach spaces.
Let us first state the  definition:

\bd Consider a Banach space $X,$ a Banach sequence space $X_d,$ and two arbitrary sequences $\etk \subset X$ and $\{e_k^*\}_{k=1}^\infty
	\in X^*.$ The pair $(\etk, \{e_k^*\}_{k=1}^\infty)$ is called an atomic decomposition of $X$ with respect
	to $X_d$, with bounds $A,B>0$, if
	\bei
	\item[(i)] $\{\la x, e_k^*\ra\}_{k=1}^\infty \in X_d$ for all $x\in X;$
	\item[(ii)]
	$ A\,||x|| \le || \{\la x, e_k^*\ra\}_{k=1}^\infty||_{ X_d}\leq B ||x||$ for all $x\in X;$
	\item[(iii)]
	$ x= \suk \la x, e_k^*\ra e_k$ for all $x\in X.$
	\eni
\ed
Note that if $\h$ is a separable Hilbert space and $X_d= \ell^2,$ the conditions
(i)+(ii) automatically imply the existence of a sequence $\etk\in \h$ such
that (iii) holds. This is well-studied in the literature on frames, see, e.g., \cite{CB}.
On the other hand, in Banach spaces the so-called {\it reconstruction property} in (iii)
does not follow from (i)+(ii) and has to be assumed separately.  Regardless whether
(i)+(ii) holds or not, the reconstruction property (iii) clearly holds if
$\etk$ is a basis for the Banach space $X$ and $\{e_k^*\}_{k=1}^\infty$ is
the dual basis. ``Genuine" atomic decompositions have been constructed, e.g., in \cite{DKST, F2, FJ,W}. Stability conditions for atomic decompositions were studied
in the paper \cite{CH}; more precisely, it was shown that if
a sequence $\ftk \subset X$ yields an atomic decomposition with
respect to a sequence $\gtk \subset X^*,$ then a sufficiently small perturbation
$\{f_k^\prime\}_{k=1}^\infty$ of $\ftk$ also yields an atomic decomposition
with respect to a certain sequence $\{g_k^\prime\}_{k=1}^\infty$. The
following result specifies how to ensure that the perturbation
condition considered in \cite{CH} is satisfied in the current setting:

\bc \label{192906a} Assume that $( \ftk,\gtk)$ is an atomic decomposition of $X$ with respect to $\ell^p$ for some $p\in[1,\infty)$, with bounds $A,B$.  Take $\epsilon_j=\epsilon 2^{-j}$ for some $0<\epsilon<B^{-1}$. Then, if \eqref{111218cn} holds, there exists a family $\{g_k^\prime\}_{k=1}^\infty\in X^*$ such that $(\{T^{\alpha(k)}\varphi\}_{k=1}^\infty , \{g_k^\prime\}_{k=1}^\infty)$ is an atomic decomposition of $X$ with respect to $\ell^p,$ with bounds $A(1+\epsilon B)^{-1}$ and $B(1-\epsilon B)^{-1}$. Moreover $\{T^{\alpha(k)}\varphi\}_{k=1}^\infty$  is a basis for $X$ if and only if $\ftk$ is a basis for $X$.
\ec
\bp  Corollary \ref{18315a} (i) implies that
if $\ftk$ satisfies the conditions stated in either Theorem \ref{18314b}, Theorem \ref{050519a} or Theorem \ref{190610d} then for every finite sequence $\{c_k\}$,
\bes
\| \sum c_k(f_k - T^{\alpha(k)}\varphi) \|_X &\leq & \sum |c_k| \|f_k - T^{\alpha(k)}\varphi \|_X \\
&\leq & (\sum |c_k|^p)^{1/p} \left( \sum \|f_k - T^{\alpha(k)}\varphi \|^q_X \right)^{1/q}\\
&\leq & \epsilon \, \| \{ c_k\} \|.
\ens The obtained inequality is a special case of the condition in
Theorem 2.3 in \cite{CH}, which immediately yields the stated conclusion.\ep

A more general way of obtaining decompositions in Banach spaces is
obtained by considering {\it Banach frames} rather than atomic
decompositions. While an atomic decomposition reconstructs an
element $x\in X$ using an infinite linear combination of the
vectors $\etk$ with coefficients $\la x, e_k^*\ra, \, k\in \mn,$ the definition
of a Banach frame ensures the existence of a certain
{\it reconstruction operator,} which maps the coefficients
$\la x, e_k^*\ra, \, k\in \mn,$ back to the vector $x\in X.$
Again based on a stability result from the paper \cite{CH} and
in a completely similar fashion as the proof of Corollary \ref{192906a},
one can prove that if a sequence $\ftk \subset X^*$ generates a Banach
frame with respect to a certain bounded operator, then an estimate
of the type \eqref{111218cn} implies that also the sequence
$\{T^{\alpha(k)}\varphi\}_{k=1}^\infty$ generates a Banach frame.
Due to the similarity with  Corollary \ref{192906a}   we will not go into details,
but just stress the fact that due to the setting of  Banach frames
the application of our theoretical results will take place in
the dual space $X^*$ for this particular case.

\section{Declarations}

\noindent Funding: None

\noindent Conflicts of interest/Competing interests: None

\noindent Availability of data and material: Not applicable

\noindent Code availability: Not applicable

\noindent Authors' contributions: All authors have contributed equally. 

 \end{document}